\documentclass[12pt]{article}

\usepackage{amsfonts}

\newtheorem{thm}{Theorem}[section]

\newtheorem{cor}[thm]{Corollary}

\newcommand{\proof
}{\par\medskip\noindent {\bf Proof.\ \ }}

\newcommand{\be}{\begin{equation}}
\newcommand{\ee}{\end{equation}}
\newcommand{\openbox}{\leavevmode
  \hbox to8pt{\hfil\vrule\vbox to6pt{\hrule width6pt\vfil\hrule}\vrule}}

\newcommand{\qed}{\hbox to5pt{ } \hfill \openbox\bigskip\medskip}

\newcommand{\Zp}{\mathbb Z _p}
\newcommand{\Zk}{\mathbb Z _k}
\newcommand{\Zm}{\mathbb Z _m}

\newcommand{\cF}{\mbox{$\cal F$}}

\newcommand{\cP}{\mbox{$\cal P$}}

\newcommand{\ve}[1]{\mathbf{#1}}

\newcommand{\N}{\mathbb N}
\newcommand{\Z}{\mathbb Z}

\title{The Erd\H{o}s-Ginzburg-Ziv  constant and progression-free subsets}
\author{G\'abor Heged\H{u}s
\\{\normalsize  \'Obuda University}
\\{\normalsize Kiscelli utca 82, Budapest, Hungary, H-1032}
\\{\normalsize hegedus.gabor@nik.uni-obuda.hu}
}

\begin{document}
\maketitle

\begin{abstract}
Ellenberg and Gijswijt gave recently a new exponential upper bound for the size of three-term arithmetic progression free sets in $({\mathbb Z _p})^n$, where $p$ is a prime. Petrov summarized their method and generalized their result to linear forms.

In this short note we use Petrov's result to give new exponential upper bounds for the   Erd\H{o}s-Ginzburg-Ziv  constant of finite Abelian groups of high rank. Our main results depend on a conjecture about Property D. 
\end{abstract}
\medskip

\noindent

\section{Introduction}


Let $A$ denote an additive finite Abelian group. 
We denote by $s(A)$ the smallest integer $\ell\in \N$ such that every sequence $S$ over $G$ of length $|S|\geq \ell$ has a zero--sum subsequence of length $|T|=exp(A)$.

Here $s(A)$ is the {\em Erd\H{o}s-Ginzburg-Ziv} constant of $A$.

Erd\H{o}s, Ginzburg and Ziv determined precisely $s(A)$ in the special case $A=\Zm$, where $m>1$ is an arbitrary integer (see \cite{EGZ}). 

Let $g(A)$ denote the smallest integer $\ell\in \N$ such that every square-free sequence $S$ over $G$ of length $|S|\geq \ell$ has a zero--sum subsequence of length $|T|=exp(A)$.

The precise value of $s(A)$  is known only for groups with rank at most two. We have
\begin{thm}
If $A=\Z_{n_1} \oplus \Z_{n_2}$, where $1\leq n_1 | n_2$, then
$s(A)=2n_1+2n_2-3$. 
\end{thm}

Let $A:=({\Zk})^n$ with $k,n\in \N$ and $k\geq 2$. The inverse problem associated with $s(A)$ asks for the structure of sequences of length $s(A)-1$ that do not have a zero-sum subsequence of length $k$. 
The standing conjecture that every group  $A:=({\Zk})^n$ satisfies the following Property D (see \cite{GG}, Conj. 7.2).

\medskip

{\bf Property D:} Every sequence $S$ over $A$ of length $|S|=s(A)-1$ that has no zero-sum subsequence of length $k$ has the form $S=T^{k-1}$ for some sequence $T$ over $A$.

\medskip

We collected the most important cases, when Property D is satisfied.
\begin{thm} \label{PropD}
The following Abelian groups has Property D:
\begin{itemize}
\item[(i)] $A=({\Zk})^n$, where $k=2$, $n\geq 1$ is arbitrary;
\item[(ii)] $A=({\Zk})^n$, where $k=3$, $n\geq 1$ is arbitrary;
\item[(iii)] $A=({\Zk})^n$, where $n=1$, $k\geq 2$ is arbitrary.
\end{itemize}
\end{thm}

Harborth proved the following inequality in \cite{H}.
\begin{thm} \label{Har}
Let $k\geq 2$, $n\geq 1$ be arbitrary integers.
Let $A:={(\Zk})^n$. Then
$$
(k-1) 2^n +1\leq s(A)\leq (k-1)k^n+1.
$$
\end{thm}

Harborth determined the Erd\H{o}s-Ginzburg-Ziv constant $s(A)$ in the following special case  in  \cite{H}.
\begin{thm} \label{Har2}
Let $a\geq 1$, $n\geq 1$ be arbitrary integers. Let $k:=2^a$. Let $A:={(\Zk})^n$. Then 
$$
s(A)=(k-1) 2^n +1.
$$
\end{thm}

Alon and Dubiner proved the following upper bound for $s(A)$ in \cite{AD}. 

\begin{thm} \label{AlDub}
Let $k\geq 2$, $n\geq 1$ be arbitrary integers. Let $A:=({\Zk})^n$. 
There exists an absolute constant $c>0$ such that for all $k$
$$
s(A)\leq (cn\mbox{log}_2 n)^nk.
$$
\end{thm}
Our main result is the following upper bound for $s(A)$, where $A:=({\Zp})^n$ and  $p>2$ is a prime.

\begin{thm} \label{upper_EGZ}
Let $p>2$ be a prime. Suppose that the group $A:=({\Zp})^n$ satisfies Property D. 

Then 
$$
s(A)\leq (p-1)p^{(1-\frac{(p-2)^2}{2 p^2ln(p)})n+1}+1.
$$
\end{thm}

We give the proof of Theorem  \ref{upper_EGZ} in Section 3.

Meshulam proved the following result in \cite{M} Corollary 1.3.

\begin{thm} \label{Mesh}
Let $n\geq 1$ be an arbitrary integer. Let $A=({{\Z}_3})^n$.  
Then
$$
s(A)=O(\frac{3^n}{n}).
$$
\end{thm}

Harborth expressed $s(A)$ in the following special case  in  \cite{H}, Hilfsatz 3.
\begin{thm} \label{Har3}
Let  $n\geq 1$ be an arbitrary integer. Let $A={{\Z}_3}^n$. Then
$$
s(A)=2g(A)-1.
$$
\end{thm}

We can confirm Alon and Dubiner's conjecture (see \cite{AD2}) using Ellenberg--Gijswijt's result about three term arithmetic progression.
\begin{cor} \label{maincor2}
Let  $n\geq 1$ be an arbitrary integer. Let $A=({{\Z}_3})^n$. 
Then 
$$
s(A)\leq  2\cdot (2.765)^{n}.
$$
\end{cor}
\proof 
Let $B\subseteq A$ be a set of vectors without three  term arithmetic progression.
It is easy to verify that 
$$
|B|\leq g(A)-1.
$$ 

The result follows from Theorem  \ref{Har3}  and  Theorem \ref{main5} . \qed

First we generalize Theorem \ref{upper_EGZ} to prime powers. 
\begin{thm} \label{ppower}
Let $p>2$ be a prime and $r\geq 1$ be an integer. Suppose that the group $A:=({\Z}_{p})^n$ satisfies Property D. Then there exists a constant $1<c(p^r)<p^r$ depending on $p^r$ such that
$$
s(({\Z}_{p^r})^n) \leq c(p^r)^n.
$$
Specially 
\be \label{EGZin}
s(({\Z}_{p^r})^n) \leq d(p)^n \frac{p^r-1}{p-1},
\ee 
for each $r\geq 1$, where $d(p)$ depends only on $p$. 
\end{thm}

The proof of Theorem \ref{ppower} appears in Section 3. 

The following result was proved in \cite{EEGKR} as Theorem 1.4.

\begin{thm} \label{EGZupper}
Let $A:={\Z}_{n_1}\bigoplus \ldots \bigoplus {\Z}_{n_r}$, where $r:=r(G)$ and $1<n_1|\ldots | n_r$. Let $c_1, \ldots ,c_r\in \N$ be integers such that for all primes $p\in \cP$ with $p| n_r$ and all $1\leq i\leq r$ we have 
$$
s({\Zp}^i)\leq c_i(p-1)+1.
$$
Then
$$
s(A)\leq \sum_{i=1}^r (c_{r+1-i}-c_{r-i})n_i -c_r +1,
$$
where $c_0=0$. In particular, if $n_1=\ldots =n_r=n$, then $s(A)\leq c_r(n-1)+1$.
\end{thm}

We use later the following easy Corollary.

\begin{cor} \label{maincor3}
Let $\cP$ denote a non-empty, finite set of odd primes and let $k\in \N$ be a product of prime powers with primes from $\cP$.

Let $m$ be a power of $2$. Suppose that for each $p\in \cP$ there exists a $1<c(p)<p$ depending only on $p$ such that
$$
s(({\Zp})^n)\leq c(p)^n(p-1)+1
$$
for each $n\geq 1$. Then there exists a $1<c(k)<k$ depending only of $k$ such that
$$
s(({\Z}_{mk})^n)\leq 2^n(m-1)k +c(k)^n (k-1) +1.
$$
\end{cor}

For the reader's convenience we give the proof of Corollary \ref{maincor3} in Section 3. But this proof is very similar to the proof of \cite{EEGKR} Corollary 4.5(3). 

Finally we  generalize Theorem \ref{upper_EGZ} to arbitrary integer modulus.


\begin{thm} \label{genEGZ} 
Let $\cP$ denote a non-empty, finite set of odd primes and let $k\in \N$ be a product of prime powers with primes from $\cP$.

Let $m$ be a power of $2$. Suppose that the groups $A_p:=({\Zp})^n$ satisfy Property D for each $p\in \cP$.
Then there exists a $1<c(k)<k$ depending only of $k$ such that
$$
s(({\Z}_{mk})^n)\leq 2^n(m-1)k +c(k)^n (k-1) +1.
$$
\end{thm}

{\bf Proof of Theorem \ref{genEGZ}:}\\

Theorem \ref{genEGZ} follows from Theorem \ref{upper_EGZ} and Corollary \ref{maincor3}. \qed




\section{Preliminaries}

\subsection{Combinatorial number theory}

Let $G$ denote a finite Abelian group. 

If $|G|>1$, then it is well--known that there exist uniquely determined integers $1<n_1 | n_2 | \ldots | n_r$ such that 
$$
G \cong {\Z}_{n_1} \oplus \ldots \oplus {\Z}_{n_r}.
$$
and $exp(G)=n_r$ is called the exponent of $G$ and $r(G)=r$ the rank of $G$.

Recall that  $G$ is a {\em $p$-group} if $exp(G)=p^k$ for a $p$ prime number and $k\in \N$, and $G$ is an {\em elementary $p$-group} if $exp(G)=p$. 

We denote by $\cF(G)$ the free abelian monoid with basis $G$. An element $S\in \cF(G)$ is called a sequence over  $G$ and we can write as:
$$
S=\prod_{g\in G} g^{\nu_g(S)}= \prod_{i=1}^l g_i,
$$  
where $\nu_g(S)\in \N$, and $l\in \N$ and $g_i\in G$. 

Here $|S|=l\in \N$ is the length, $\sigma(S)=\sum_{i=1}^l g_i\in G$ is the sum and $supp(S)=\{g\in G:~ \nu_g(S)>0\}$ is the support of $S$.  Moreover,  $\nu_g(S)$ is called the multiplicity of $g$ in $S$. 

Recall that a sequence is called a {\em zero-sum sequence} if $\sigma(S)=0$, it is called squarefree if $\nu_g(S)\leq 1$ for each $g\in G$. As usual, a sequence $T$ is called a subsequence  of $S$ if $T$ divides $S$ in $\cF(G)$. 

We use later the following result (see \cite{CDGGS} Proposition 3.1).
\begin{thm} \label{exp_upper}
Let $G$ be a finite Abelian group and let $H\leq G$ be a subgroup such that $exp(G)=exp(H)exp(G/H)$. Then
$$
s(G)\leq exp(G/H)(s(H)-1)+s(G/H)
$$
\end{thm}

\subsection{Upper bounds for the space of monomials}

Let $n,m\geq 1$, $D\geq 2$ be fixed integers. Let 
$$
L_{n,D}:=\mbox{span}(\{x^{\ve\alpha}=x_1^{\alpha_1}\ldots x_n^{\alpha_n}:~{\alpha_i\leq D-1} \mbox{ for each } 1\leq i\leq n\} ).
$$
Let 
$$
L_{n,D,k}:=\mbox{span}(\{x^{\ve\alpha}\in L_{n,D}  :~ deg(x^{\ve\alpha})\leq k \}).
$$
\begin{thm} \label{upper33}
Let 
$$
c:=1-\frac{(m-2)^2}{2 m^2ln(D)}.
$$
Then 
$$
dim(L_{n,D,\frac{n(D-1)}{m}})\leq D^{cn}.
$$
\end{thm}
The proof of Theorem \ref{upper33} is based on Hoeffding inequality and a very slight generalization of the proof of \cite{G} Lemma 1. 

\subsection{Progression-free sets and Petrov's result}

Ellenberg and Gijswijt achieved the following breakthrough very recently in \cite{EG} on the upper bounds of progression free subsets in $({\Zp})^n$, where $p$ is a prime.

\begin{thm} \label{main55}
Let $p$ be a fixed prime.
Let $(a_1,a_2 ,a_3)\in ({\Zp})^3$ be  a fixed vector. Suppose that 
$p$ divides $\sum_i a_i$ and  $a_3\not\equiv 0 \pmod p$.. 

Let $\cF\subseteq ({\Zp})^n$ be an arbitrary  subset satisfying the following property:   if $\ve b_1, \ve b_2, \ve b_3\in {\cF}$ are arbitrary vectors such that 
there exist $1\leq i<j\leq 3$ with $\ve b_i\ne \ve b_j$, then
$$
a_1\ve b_1+ a_2\ve b_2+ a_3\ve b_3\ne \ve 0.
$$

Then 
$$
|\cF|\leq 3 \cdot dim(L_{n,p,\frac{n(p-1)}{3}}).
$$
\end{thm}

Later Petrov proved the following generalization of the  Theorem  \ref{main55} in \cite{P}. 
\begin{thm} \label{main4}
Let $p$ be a prime, $m\geq 1$ be an integer.
Let $(a_1,\ldots ,a_m)\in ({\Zp})^m$ be  a fixed vector. Suppose that 
$p$ divides $\sum_i a_i$ and  $a_m\not\equiv 0 \pmod p$.

Let $\cF\subseteq ({\Zp})^n$ be an arbitrary subset satisfying the following property:  if $\ve b_1, \ldots, \ve b_m\in {\cF}$ are arbitrary vectors such that 
there exist $1\leq i<j\leq m$ with $\ve b_i\ne \ve b_j$, then
$$
\sum_{i=1}^m a_i\ve b_i \ne \ve 0.
$$

Then 
$$
|\cF|\leq m \cdot dim(L_{n,p,\frac{n(p-1)}{m}}).
$$
\end{thm}

The following Corollary is a direct consequence of Theorem \ref{upper33}.
\begin{cor} \label{main5}
Let $p$ be a prime, $m\geq 1$ be an integer.
Let $(a_1,\ldots ,a_m)\in ({\Zp})^m$ be  a fixed vector. Suppose that 
$p$ divides $\sum_i a_i$. 

Let $\cF\subseteq ({\Zp})^n$ be an arbitrary subset satisfying the following property:  if $\ve b_1, \ldots, \ve b_m\in {\cF}$ are arbitrary vectors such that 
there exist $1\leq i<j\leq m$ with $\ve b_i\ne \ve b_j$, then
$$
\sum_{i=1}^m a_i\ve b_i \ne \ve 0.
$$

Then 
$$
|\cF|\leq mp^{(1-\frac{(m-2)^2}{2 m^2ln(p)})n}.
$$
\end{cor}
\proof 
This follows from Theorem \ref{main4} and Theorem \ref{upper33}. \qed
Finally we use in the proof of our main results the following clear consequence. 

\begin{cor} \label{main7}
Let $p$ be a prime, $r\geq 2$ be an integer. Suppose that 
$p$ divides $r$. 

Let $\cF\subseteq ({\Zp})^n$ be an arbitrary subset satisfying the following property:  if $\ve b_1, \ldots, \ve b_r\in {\cF}$ are arbitrary vectors such that 
there exist $1\leq i<j\leq r$ with $\ve b_i\ne \ve b_j$, then
$$
\ve b_1+ \ldots +\ve b_r \ne \ve 0.
$$

Then 
$$
|\cF|\leq rp^{(1-\frac{(r-2)^2}{2 r^2ln(p)})n}.
$$
\end{cor}
\proof 
Let $m:=r$ and $a_i=1$ for each $1\leq i\leq r$. Then we
can apply Corollary \ref{main5}.
\qed

\section{Proofs}

{\bf Proof of Theorem \ref{upper_EGZ}:}\\

Let $S$ be a sequence of length $s(A)-1$ that does not have a zero-sum subsequence of length $p$. Then  
it follows from Property D that $S$ has the form $S=T^{p-1}$ for some sequence $T$ over $A$. Clearly $T$ does not have a zero-sum subsequence of length $p$. Hence $T$ is a square-free sequence, and there exists a 
subset  $\cF\subseteq {\Zp}^n$ such that $\cF=supp(T)$. Hence  if $(\ve a_1, \ldots, \ve a_{r})\in {\cF}^{p}$ is an arbitrary vector such that 
there exist $1\leq i<j\leq m$ with $\ve a_i\ne \ve a_j$, then 
$$
\ve a_1+ \ldots +\ve a_p \ne \ve 0.
$$

Here we used that $S=T^{p-1}$ is a sequence of length $s(A)-1$ that does not have a zero-sum subsequence of length $p$.

We can apply Corollary \ref{main7} with the choice $r:=p$ and we get
$$
|\cF|\leq p\cdot p^{(1-\frac{(p-2)^2}{2 p^2ln(p)})n}.
$$
Since  $S=T^{p-1}$, $\cF=supp(T)$ and $S$ is a sequence of length $s(A)-1$, hence we get our result.  
\qed

{\bf Proof of Theorem \ref{ppower}:}\\

We can prove by induction on the exponent $r$.

If $r=1$, then Theorem \ref{upper_EGZ} gives the result.

Suppose that our inequality (\ref{EGZin}) is true for  a  fixed $r$. We will prove that (\ref{EGZin}) is true for $r+1$. 
Namely it follows from Theorem \ref{upper_EGZ} and the inductional hypothesis that 
$$
s(({\Z}_{p^{r+1}})^n) \leq p\cdot s(({\Z}_{p^r})^n) + s(({\Z}_{p})^n)\leq 
$$
$$
\leq d(p)^n \frac{p^r-1}{p-1} + d(p)^n=   d(p)^n \Big(  \frac{p^r-1}{p-1}  +1\Big)=
$$
$$
d(p)^n \frac{p^{r+1}-1}{p-1}.
$$
\qed

{\bf Proof of Corollary \ref{maincor3}:}\\

Let $A:=({\Z}_{mk})^n$ and $H:=kA \equiv ({\Z}_{m})^n$. Then clearly $A/H \equiv ({\Z}_{k})^n$. 
Define
$$
c(k):=\mbox{max}\{c(p):~ p\in \cP \}.
$$
Clearly $1<c(k)<k$, since $1<c(p)<p$ for each $p\in \cP$.

But then 
$$
s({\Zp})\leq \ldots \leq s(({\Zp}^n)) \leq c(p)^n(p-1)+1\leq c(k)^n(p-1)+1
$$
for each $p\in \cP$. It follows from Theorem \ref{EGZupper} that 
$$
s(A/H)=s(({\Z}_{k})^n)\leq c(k)^n(k-1)+1.
$$

On the other hand it follows from Theorem \ref{Har2} that
$$
s(H)=s(({\Z}_{m})^n) \leq 2^n(m-1)+1.
$$

Finally Theorem \ref{exp_upper} gives that
$$
s(A)\leq exp(A/H)(s(H)-1)+s(A/H)\leq 2^n(m-1)k + c(k)^n(k-1)+1.
$$
\qed

{\bf Acknowledgements.} 
I am indebted to  Lajos R\'onyai for his useful remarks.

\end{document}